\documentclass[10pt,twocolumn,twoside]{IEEEtran}
\usepackage{etex}

\IEEEoverridecommandlockouts

\usepackage{amsmath,mathrsfs, enumerate}
\usepackage{amsthm}
\usepackage{amssymb}
\usepackage{mathtools}
\usepackage{pstricks,pst-plot,psfrag}
\usepackage[all]{xy}
\usepackage{graphicx,xspace,bm}
\usepackage{color}
\usepackage{hyperref}

\usepackage{subfigure}

\usepackage{algorithm}
\usepackage[noend]{algpseudocode}

\usepackage{epstopdf}         
\graphicspath{{./epsfiles/}} 
\usepackage{stackengine}
 
\allowdisplaybreaks
 
\newtheorem{theorem}{Theorem}[section]

\newtheorem{problem}[theorem]{Problem}

\newtheorem{corollary}[theorem]{Corollary}
\newtheorem{proposition}[theorem]{Proposition}

\newcommand{\realnonnegative}{{\mathbb{R}}_{\ge 0}}
\newcommand{\integers}{\mathbb{Z}}
\newcommand{\integerspositive}{\mathbb{Z}_{\geq 1}}

\newcommand{\oprocendsymbol}{\hbox{$\square$}}
\newcommand{\oprocend}{\relax\ifmmode\else\unskip\hfill\fi\oprocendsymbol}

\DeclareMathOperator*{\argmin}{arg\,min}

\newcommand{\longthmtitle}[1]{\mbox{}\textup{\bf (#1):}}

\newcommand{\ra}[1]{\renewcommand{\arraystretch}{#1}}

\newcommand{\Zprocess}[1]{\{Z_{#1}\}_{n=1}^{\infty}}

\newcommand{\Nsum}{\sum_{i=1}^N}
\newcommand{\tsum}{\sum_{t=1}^n}
\newcommand{\jsum}{\sum_{j\in\N_i'}}

\newcommand{\bR}{\bar{R}}
\newcommand{\bB}{\bar{B}}
\newcommand{\bT}{\bar{T}}

\newcommand{\budget}{\mathcal{B}}
\newcommand{\F}{\mathcal{F}}
\newcommand{\N}{\mathcal{N}}
\newcommand{\G}{{\mathcal{G}}}
\newcommand{\E}{{\mathcal{E}}}
\newcommand{\X}{{\mathcal{X}}}

\newcommand{\avginf}{\tilde{I}_n}
\newcommand{\suscept}{\tilde{U}_n}
\newcommand{\exposure}{\tilde{S}_n}

\newcommand{\stratUn}{\text{(i)}}
\newcommand{\stratSn}{\text{(ii)}}
\newcommand{\stratgf}{\text{(iii)}}
\newcommand{\stratheur}{\text{(iv)}}
\newcommand{\stratunif}{\text{(v)}}

\DeclarePairedDelimiter{\ceil}{\lceil}{\rceil}

\newcommand{\filt}{\{\F_n\}_{n=1}^{\infty}}

\usepackage{flushend}

\begin{document}

\title{Curing Epidemics on Networks using a Polya Contagion Model}

\author{Mikhail Hayhoe$^1$ \qquad Fady Alajaji$^2$ \qquad Bahman Gharesifard$^2$
\thanks{
$^1$Department of Electrical and Systems Engineering at the University of Pennsylvania, \texttt{mhayhoe@seas.upenn.edu}.

$^2$Department of Mathematics \& Statistics at Queen's University, Kingston, Ontario, Canada,~\texttt{\{fady, bahman\}@mast.queensu.ca}.

This work was partially supported by the Natural Sciences and Engineering Research Council of Canada. Parts of this work were submitted for presentation at the 2018 American Control Conference.}
}
\maketitle

\begin{abstract}
We study the curing of epidemics of a network contagion, which is modelled using a variation of the classical Polya urn process that takes into account spatial infection among neighbouring nodes. We introduce several quantities for measuring the overall infection in the network and use them to formulate an optimal control problem for minimizing the average infection rate using limited curing resources. We prove the feasibility of this problem under high curing budgets by deriving conservative lower bounds on the amount of curing per node that turns our measures of network infection into supermartingales. We also provide a provably convergent gradient descent algorithm to find the allocation of curing under limited budgets. Motivated by the fact that this strategy is computationally expensive, we design a suit of heuristic methods that are locally implementable and nearly as effective. Extensive simulations run on large-scale networks demonstrate the effectiveness of our proposed strategies.
\end{abstract}
\vspace{0.3cm}

\section{Introduction}\label{sec:intro}

In this paper we examine the problem of curing an epidemic using a network contagion model adapted from the Polya process~\cite{MH-FA-BG:17, MH-FA-BG:17-2}. Here an epidemic can represent a disease~\cite{LK-MA-KD-SK-AO:14}, a computer virus~\cite{MG-WG-DT:03}, the spread of an innovation, rumour or idea~\cite{ER:03}, or the dynamics of competing opinions in a social network~\cite{EA-LAA:05}.


Epidemics on networks have been intensively studied in recent years, see~\cite{PVM-JO-RK:09, CN-VMP-GJP:16} and references therein and thereafter. Our model is similar to the well-known susceptible-infected-susceptible (SIS) compartmental infection model~\cite{DE-JK:10}, in the sense that initially, all nodes may be healthy or infected and as the epidemic spreads, nodes that are infected can be cured to become healthy, but any healthy node may become infected at any time, regardless of whether they have been cured previously. However, the dynamics of the traditional SIS model tend to be complicated, and often deterministic approximation methods are employed to simplify the analysis~\cite{CN-VMP-GJP:16}. In contrast to the SIS model, our model is motivated by the classical Polya contagion process, which evolves by sampling from an urn containing a finite number of red and black balls~\cite{GP-FE:23,GP-FE:28,GP:30}. In the network Polya contagion model, introduced in~\cite{MH-FA-BG:17}, each node of the underlying network is equipped with an individual urn; however, instead of sampling from these urns when generating its contagion process, each node has a ``super urn'', created by combining the contents of its own urn with those of its neighbours' urns. This adaptation captures the concept of spatial infection, since having infected neighbours increases the chance that an individual is infected in the future. The stochastic properties of this model were examined in~\cite{MH-FA-BG:17, MH-FA-BG:17-2}. In this work, we study the problem of controlling the contagion spread under this model.

More specifically, we propose various natural ways to measure the total infection in the network Polya contagion model, and examine conditions under which these measures have limits as time grows without bound. Using these measures, we pose an optimal control problem within the context of the network Polya contagion model. As our first contribution, we characterize lower bounds on the allocation of curing to individual nodes which turn these infection measures into supermartingales. Our result hence provides a conservative strategy for curing network epidemics. We next focus on realistic scenarios, where the curing budget is constrained. As our next contribution, we prove that the constrained gradient flow method is convergent for this problem and hence can be employed to find near-optimal strategies under a fixed curing budget at each time step. In spite of its effectiveness, as we demonstrate, the gradient flow strategy is computationally expensive and is only implementable in a centralized manner. These shortcomings motivate us to look into alternative strategies, which take advantage of notions of node centrality of the underlying network along with the composition of super urns at each time step. These strategies are less expensive computationally and can be adopted for implementation in a decentralized manner. Through extensive simulation results, we show that our proposed heuristic strategies perform well in curing epidemics.

%

The rest of the paper is organized as follows. Section~\ref{sec:prelim} outlines some mathematical preliminaries that will be used throughout the paper. Section~\ref{sec:model} contains the description of our network contagion process and the problem statement. Section~\ref{sec:control_analytical} discusses analytical results pertaining to the control of epidemics using this model, and Section~\ref{sec:curing} outlines strategies used to solve the problems posed. Section~\ref{sec:simulations} includes several simulation results. Finally, Section~\ref{sec:conclusion} sumamrizes our contributions and proposes avenues for future work.

\section{Preliminaries}\label{sec:prelim}

For a sequence $v_i = (v_{i,1},\ldots,v_{i,n})$, we use the notation $v_{i,s}^t$ with $1 \leq s < t \leq n$ to denote the vector $(v_{i,s},v_{i,s+1},\ldots,v_{i,t})$, with $ v_{i,0}^t = v_i^t $. Our technical results rely on notions from stochastic processes, some of which we recall here. Throughout, we assume that the reader is familiar with basic notions of probability theory. 

Let $ (\Omega,\F,P) $ be a probability space, and consider the stochastic process $\Zprocess{n}$, where each
$ Z_n $ is a random variable on $\Omega$. We often refer to the indices of the process as ``time'' indices. We recall that the process $\Zprocess{n}$ is \emph{stationary} if for any $n\in \integers_{\geq 1}$, its $n$-fold joint probability distribution (i.e., the distribution of $(Z_1,\ldots,Z_n)$) is invariant to time shifts. Further, $\Zprocess{n}$ is \emph{exchangeable} if for any $n\in \integers_{\geq 1}$, its $n$-fold joint distribution is invariant to permutations of the indices $1,\ldots,n$. It directly follows from the definitions that an exchangeable process is stationary. Lastly, the process $\Zprocess{n}$ is called a \emph{martingale} (resp. \emph{supermartingale}, \emph{submartingale}) with respect to the filtration $ \{\F_n\}_{n=1}^{\infty} $ if $E[|Z_n|] < \infty$ and $E[Z_{n+1} | \F_{n}] = Z_{n}$ almost surely (resp. less than or equal to, greater than or equal to), for all $n$. If the inequality is strict, we call the process a \emph{strict supermartingale} or \emph{strict submartingale}. Doob's martingale convergence theorem~\cite{RA-CD:00} can then be used to show that $\Zprocess{n}$ will have a limit as $ n $ grows without bound. Precise definitions of all notions, including that of {\em ergodicity}, can be found in standard texts (e.g.,~\cite{RA-CD:00,GG-DS:01}).

\begin{figure}[!ht]
\centering 
	\includegraphics[width=\linewidth]{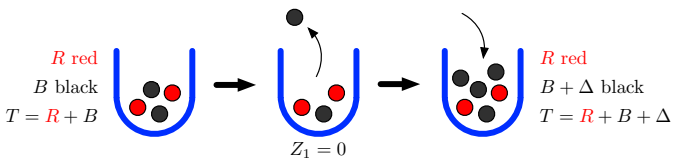}
	\caption{Illustration of the first draw for a classical Polya process. We drew a black ball and hence $ Z_1 = 0 $. Here $ R = 2 $, $ B = 2 $, and $ \Delta = 2 $.}
	\label{fig:urn_draw}
\end{figure}

We now recall the classical version of the Polya contagion process~\cite{GP-FE:23,GP:30}. Consider an urn with $R \in \integers_{> 0} $ red balls and $B \in \integers_{> 0}$ black balls. We denote the total number of balls by $ T $, i.e., $T = R + B$. At each time step, a ball is drawn from the urn. The ball is then returned along with $\Delta\ > 0$ balls of the same color. To describe this process, we use a replacement matrix $ M_R $:
\[
	M_R = 	\begin{bmatrix}
				\Delta & 0 \\
				0 & \Delta
			\end{bmatrix},
\]
where $ [M_r]_{1,1} = \Delta, [M_r]_{1,2} = 0 $ means we add $ \Delta $ red balls and $ 0 $ black balls when a red ball is drawn, and similarly $ [M_r]_{2,1} = 0, [M_r]_{2,2} = \Delta $ means we add $ 0 $ red balls and $ \Delta $ black balls when a black ball is drawn. and We use an indicator $Z_n$ to denote the color of ball in the $n$th draw (see Figure~\ref{fig:urn_draw}):
\begin{align*}
	Z_n = \begin{cases}
		1 &\text{if the $n$th draw is red}\\
		0 &\text{if the $n$th draw is black.}
	\end{cases}
\end{align*}

Let $U_n$ denote the proportion of red balls in the urn after the $n$th draw. Then
\begin{align*}
	U_n &:= \frac{R + \Delta\tsum Z_{t}}{T + n\Delta}\cr
		&= \frac{\rho_c + \delta_c\tsum Z_{t}}{1+n\delta_c}
\end{align*}
where $\rho_c = \frac{R}{T}$ is the initial proportion of red balls in the urn and $\delta_c = \frac{\Delta}{T}$ is a correlation parameter. Since we draw balls from this urn at each time step, the conditional probability of drawing a red ball at time $n$, given $Z^{n-1}=(Z_1,\ldots,Z_{n-1})$, is given by
\begin{align*}
	P(Z_n = 1 \ | \ Z^{n-1} ) &= \frac{R + \Delta\sum_{t=1}^{n-1} Z_{t}}{T + (n-1)\Delta} \\
		&= U_{n-1}.
\end{align*}
It can be easily shown that $\{U_n\}_{n=1}^{\infty}$ is a martingale~\cite{WF:71}.
The process $\{Z_n\}_{n=1}^{\infty}$, whose $n$-fold joint distribution can be determined in closed form in terms of the parameters $ \rho_c $ and $ \delta_c $, is also exchangeable (hence stationary) and non-ergodic with both $U_n$ and the process sample average $\frac{1}{n}\sum_{i=1}^n Z_i$ converging almost surely as $n \rightarrow \infty$ to a random variable governed by the Beta distribution with parameters $\frac{\rho_c}{\delta_c}$ and $\frac{1-\rho_c}{\delta_c}$~\cite{WF:71, FA-TF:94}. The classical Polya process has been applied in many different contexts, including the modelling of communication channels with memory~\cite{FA-TF:94}, image segmentation~\cite{AB-PB-FA:99}, as well as in biology, statistics and other areas (see~\cite{RP:07}).

\section{Model Description and Problem Statement}\label{sec:model}

\subsection{Network Polya Contagion Process}\label{subsec:process}
In this section, we briefly recall the Polya network contagion process introduced in~\cite{MH-FA-BG:17, MH-FA-BG:17-2}. Consider an undirected graph $\G = (V,\E) $, where $V=\{1,\ldots, N\} $ is the set of $ N \in \integers_{\geq 1} $ nodes and $ \E \subset V\times V $ is the set of edges. We assume that $ \G $ is connected, i.e. there is a path between any two nodes in $ \G $. We use $\N_i$ to denote the set of nodes that are neighbors to node $i$, that is $\N_i = \{v \in V : (i,v) \in \E\}$, and $\N_i' = \{i\} \cup \N_i$. Each node $i\in V$ is equipped with an urn,  initially with $R_i \in \integers_{> 0}$ red balls and $B_i \in \integers_{> 0}$ black balls (we do not let $R_i = 0$ or $B_i=0$ to avoid any degenerate cases). We let $T_i = R_i + B_i$ be the total number of balls in the $ i $th urn, $i \in \{1,\ldots,N\}$. We use $Z_{i,n}$ as an indicator for the ball drawn for node $i$ at time $n$:
\begin{align*}
	Z_{i,n} = \begin{cases}
		1 &\text{if the $n$th draw for node $i$ is red}\\
		0 &\text{if the $n$th draw for node $i$ is black.}
	\end{cases}
\end{align*}
Thus we define the network contagion process as $ \{Z_n\}_{n=1}^{\infty} $, where $ Z_n = (Z_{i,n},\ldots,Z_{N,n}) $.  However, instead of drawing solely from its own urn, each node has a ``super urn'' created by combining all the balls in its own urn with the balls in its neighbours' urns; see Figure~\ref{fig:super_urn}. This allows the spatial relationships between nodes to influence their state. This means that $Z_{i,n}$ is the indicator for a ball drawn from node $i$'s super urn, and not its individual urn. Hence, the super urn of node $ i $ initially has $\bR_i = \sum_{j \in \N_i'}R_j$ red balls, $\bB_i = \sum_{j \in \N_i'}B_j$ black balls, and $\bT_i = \sum_{j \in \N_i'}T_j$ balls in total, and the network-wide initial proportion of red balls is $ \rho = \frac{\Nsum R_i}{\Nsum T_i} $.
{
\psfrag{1}[rr][rr]{{\footnotesize $1$}}
\psfrag{2}[rr][rr]{{\footnotesize $2$}}
\psfrag{3}[rr][rr]{{\footnotesize $3$}}
\psfrag{4}[rr][rr]{{\footnotesize $4$}}
\psfrag{5}[rr][rr]{{\footnotesize $5$}}
\psfrag{6}[rr][rr]{{\footnotesize $6$}}
\psfrag{7}[rr][rr]{{\footnotesize $7$}}
\psfrag{Node 1's super urn}[rr][rr]{{\footnotesize Node 1's super urn}}
\begin{figure}[!ht]
\centering 
	\includegraphics[width=0.8\linewidth]{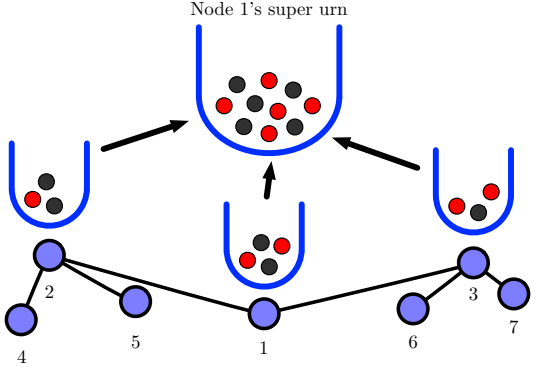}
	\caption{Illustration of a super urn in a network.}
	\label{fig:super_urn}
\end{figure}
}

We allow the number of added balls to vary based on the colour drawn, the time index, and the node for which it was drawn; hence, the replacement matrix for node $ i $ at time $ t $ is
\[
	M_{R,i}(t) = 	\begin{bmatrix}
				\Delta_{r,i}(t) & 0 \\
				0 & \Delta_{b,i}(t)
			\end{bmatrix}.
\]
We assume that $\Delta_{r,i}(t) \geq 0 $ and $ \Delta_{b,i}(t) \geq 0$ for all $t \in \integerspositive $, and that there exists $i\in V$ and $ t $ such that $\Delta_{r,i}(t) + \Delta_{b,i}(t) \neq 0$; otherwise we are simply sampling with replacement. We assume throughout that $\Delta_{r,i}(t) \geq 0, \Delta_{b,i}(t) \geq 0$, for all $t \in \integerspositive $ and that there exists $i\in V$ and $ t $ such that $\Delta_{r,i}(t) + \Delta_{b,i}(t) \neq 0$; otherwise we are simply sampling with replacement. In the context of epidemics, the red and black balls in an urn, respectively, represent ``infection'' and ``healthiness''. We refer the interested reader to~\cite{MH-FA-BG:17-2} for a complete description of this relationship.

To express the proportion of red balls in the individual urns of the nodes, we define the random vector $U_n=(U_{1,n},\ldots, U_{N,n})$, where $U_{i,n}$ is the proportion of red balls in node $i$'s urn after the $n$th draw, $ i \in V $. For node $i$, 
\begin{align*}
	U_{i,n} &= \frac{R_i + \tsum \Delta_{r,i}(t)Z_{i,t}}{X_{i,n}},
\end{align*}
where 
\begin{align}\label{eq:X_n}
X_{i,n}=T_i + \sum_{t=1}^{n} \Delta_{r,i}(t)Z_{i,t} + \Delta_{b,i}(t)(1-Z_{i,t})
\end{align}
is the total number of balls in node $i$'s urn after the $n$th draw, and the numerator represents the total number of red balls  in the same urn. We now define the random vector $ S_n = (S_{1,n},\ldots,S_{N,n}) $ as the proportion of red balls in the super urns of the nodes after the $ n $th draw, so that $ S_{i,n} $ is the proportion of red balls in node $ i $'s super urn after $ n $ draws. Hence, for node $ i $,
\begin{align}\label{eq:S_n}
	S_{i,n} 
	&= \frac{\bR_i + \sum_{t=1}^{n}\bar{Z}_{r,i,t}}{\bar{X}_{i,n}}  \cr
	&= \frac{\sum_{j \in \N_i^{'}} U_{j,n}X_{j,n}}{\bar{X}_{i,n}},
\end{align}
where
\begin{align*}
	\bar{Z}_{r,i,n} &= \jsum\Delta_{r,j}(n)Z_{j,n}, \\
	\bar{X}_{i,n} &= \bar{T}_i + \sum_{t=1}^n (\bar{Z}_{r,i,t} + \bar{Z}_{b,i,t}) = \jsum X_{j,n}, \\
	\bar{Z}_{b,i,n} &= \jsum\Delta_{b,j}(n)(1 - Z_{j,n}).
\end{align*}
Note that $ S_{i,0} = \bR_i/\bT_i $. In fact, $ S_{i,n} $ is a function of the random draw variables of the network, and in particular of $ \{Z_j^{n}\}_{j \in \N_i'} $, but for ease of notation, when the arguments are clear, we write $ S_{i,n}(Z_1^{n},\cdots,Z_N^{n}) = S_{i,n} $. Then the conditional probability of drawing a red ball from the super urn of node $i$ at time $n$ given the complete network history, i.e. given all the past $n-1$ draw variables for each node in the network $\{Z_j^{n-1}\}_{j=1}^N=\{(Z_{1,1},\ldots,Z_{1,n-1}),\ldots,(Z_{N,1},\ldots,Z_{N,n-1})\}$, satisfies  
\begin{align}\label{eq:cond_eq}
	P&\left(Z_{i,n} = 1 | \{Z_j^{n-1}\}_{j=1}^N\right) \cr
	&= \frac{\bR_i + \sum_{t=1}^{n-1} \bar{Z}_{r,i,n}}{\bar{X}_{i,n}} \nonumber \\
	&= S_{i,n-1}.
\end{align}
That is, the conditional probability of drawing a red ball for node $ i $ given the entire past $ \{Z_j^{n-1}\}_{j=1}^N $ is the proportion of red balls in its super urn, $ S_{i,n-1} $. Since these random variables fully describe the evolution of the process, we say $ \filt $ is the natural filtration on $ \{Z_i^{n-1}\}_{i=1}^N $ and by extension $ \{U_{i,n}\}_{n=1}^{\infty} $ and $ \{S_{i,n}\}_{n=1}^{\infty} $, for all $ i \in V $. Thus note that in~\eqref{eq:cond_eq} we could have instead conditioned on $ \F_{n-1} $.

Using the conditional probability established above, we next determine the $n$-fold joint probability of the entire network $\G$. Let $ a_i^{n} \in \{0,1\}^n$, where $i \in \{1,\ldots,N\}$. We have
\begin{align}
	&P_{\G}^{(n)}(a_1^n,\cdots,a_N^n) \cr
	&:= P\left(\{Z_i^n = a_i^n\}_{i=1}^N\right) \cr
	&= \prod_{t=1}^nP\left(\{Z_{i,t}=a_{i,t}\}_{i=1}^N\ | \ \{Z_i^{t-1} = a_i^{t-1}\}_{i=1}^N\right) \cr
	&= \prod_{t=1}^n\prod_{i=1}^N \Big(S_{i,t-1}\Big)^{a_{i,t}}\Big(1 - S_{i,t-1}\Big)^{1-a_{i,t}},
	\label{eq:joint_eq}
\end{align}
where $S_{i,t} = S_{i,t}(a_1^t,\cdots,a_N^t) $ is defined in \eqref{eq:cond_eq}. The study of the asymptotic behaviour of each node's contagion process $ \Zprocess{i,n}, i \in V $ is established in~\cite{MH-FA-BG:17, MH-FA-BG:17-2}. Our objective in this work is to demonstrate the implications of these results in designing curing strategies. With the above explicit joint distribution, it is possible to determine the distributions of each node's process. More specifically, using~\eqref{eq:joint_eq}, the $n$-fold distribution of node $i$'s process at time $t \geq n$ is
\begin{align*}
P_{i,t}^{(n)}&(a_{i,t-n+1},\cdots,a_{i,t}):= \sum_{\substack{a_i^{t-n} \in \{0,1\}^{t-n} \\ a_j^t \in \{0,1\}^t, j\ne i}} P_{\G}^{(n)}(a_1^t,\cdots,a_N^t).
\end{align*} 
It can be shown that the draw contagion process $\Zprocess{i,n}$ of each individual node $i$  is not stationary in general (and hence not exchangeable) ~\cite{MH-FA-BG:17, MH-FA-BG:17-2}. Thus the entire network contagion process $ \Zprocess{n} $ is not stationary.

In order to measure the spread of contagion in the network at any given time, we wish to see how likely it is, on average, for a node to be infected at that instant. We thus define the {\em average infection rate} in the network at time $n$ as
\[
\avginf := \frac{1}{N} \Nsum P(Z_{i,n}=1) = \frac{1}{N} \Nsum P_{i,n}^{(1)}(1).
\]
Note that $ \avginf $ is a function of the network topology $ (V,\E) $, the initial placement of balls $ R_i $ and $ B_i $, the draw processes $ \{Z_{i,t}\}_{t=1}^n $, and number of balls added $ \{\Delta_{r,i}(t)\}_{t=1}^n $ and $ \{\Delta_{b,i}(t)\}_{t=1}^n $ for each node $ i \in V $. Unfortunately for an arbitrary network, the above quantity does not yield an exact analytical formula (except in the simple case of complete networks). As such, in general it is hard to mathematically analyze the asymptotic behavior of $ \avginf $, which we wish to minimize when attempting to cure an epidemic. Instead we examine the asymptotic stochastic behavior of two closely related variables given by the average individual proportion of red balls at time $n$, namely
\[
	\suscept  := \frac{1}{N} \Nsum U_{i,n}, 
\]
which we call the {\em network susceptibility}, and the average neighborhood proportion of red balls at time $ n $,
\[
	\exposure := \frac{1}{N} \Nsum S_{i,n},
\]
which we call the {\em network exposure}.

With the model in hand, we turn to the exploration of a curing problem. Our objective is to control the average infection rate $ \avginf $ as $n$ grows without bound; but when seeking analytic results, it might be more amenable to observe the asymptotic behavior of the network exposure $ \exposure $. These quantities are closely related; through~\eqref{eq:S_n} we see that if $ U_{i,n} $ increases then this node-specific value causes $ S_{j,n} $ to increase for every neighbour $ j $ of node $ i $, and hence by~\eqref{eq:cond_eq} their conditional probabilities of drawing red balls increase. More specifically,
\begin{align}\label{rmk:measures}
	\uparrow U_{i,n} &\xRightarrow{\eqref{eq:S_n}} \quad \uparrow S_{j,n} \text{ for all } j \in \N_i' \cr
	&\xRightarrow{\eqref{eq:cond_eq}} \quad \uparrow P\left(Z_{i,n+1} = 1 | \{Z_j^{n}\}_{j=1}^N\right) \text{ for all } j \in \N_i'. \quad
\end{align}
Thus if $ \suscept $ is high, then this average measure of individual nodes implies that the conditional probability of a node being infected is higher on average. Hence $ \suscept $ can be understood as the \textit{average node prevalence} of infection. The effect of the network exposure here is more direct, since~\eqref{eq:cond_eq} shows that $ \exposure $ is in fact the network-wide average of the conditional probabilities of infection, which is a quantity that is intimately related to the state of infection in the neighbourhood of node $ i $. Thus $ \exposure $ represents the \textit{average neighbourhood prevalence} of infection. Note that similarly to $ \avginf $, both $ \suscept $ and $ \exposure $ are functions of the network variables.

\subsection{Establishing a Control Problem}
The quantities $ \{\Delta_{b,i}(n)\}_{n=1}^{\infty} $, which denote the net number of ``healthy'' balls added to node $i$'s urn after each draw, can play the role of ``healing or curing  parameters''. Our objective is to show that when these parameters are appropriately selected, one can steer the average infection rate towards a desirable level; the selection of curing parameters is, however, subject to an allowable budget on the maximal number of healthy balls that can be added in the network. Let us state this problem formally.


\begin{problem}\longthmtitle{Average Infection Rate Budget Constraint}\label{prob:budget}
Minimize the limiting average infection rate $\tilde{I}_t$ subject to a budget $\budget$ on the total healing at each time step:%
\[
\min_{\substack{\Nsum \Delta_{b,i}(t) \leq \budget \\ \forall t}} \limsup_{t\to\infty}\tilde{I}_t
\]
\end{problem}

%
%


Such optimal curing problems have been studied in many different contexts~\cite{ERL-SM:16, CN-VMP-GJP:17}. For our model, the solution to Problem~\ref{prob:budget} would be an infinite horizon optimal control policy that would yield the best possible level of epidemic elimination, given the initial data. Finding such a policy in general appears to be difficult. Nevertheless, as we demonstrate in the upcoming sections, one can obtain interesting analytical results regarding the feasibility of this problem, and design algorithmic strategies to curtail the average infection rate.

\section{Controlling Epidemics: Analytical Results}\label{sec:control_analytical}

In order to determine when Problem~\ref{prob:budget} makes sense, we wish to examine when a limit exists. As stated earlier, working with $ \avginf $ can be difficult, and so we instead focus on the related measures of the network susceptibility $ \suscept $ and network exposure $ \exposure $. Our next results will show how one can force these measures to form supermartingales by appropriately selecting the curing policies $ \{\Delta_{b,i}(n)\}_{n=1}^{\infty} $, for all $ i \in V $. In conjunction with Doob's martingale convergence theorem~\cite{RD:96}, these results show that $ \{U_{i,n}\}_{n=1}^{\infty} $, $ \{S_{i,n}\}_{n=1}^{\infty} $, and hence both $ \{\suscept\}_{n=1}^{\infty} $ and $ \{\exposure\}_{n=1}^{\infty} $, have limits. While the results presented herein do not obey the per-step budget constraint on the curing, these results in conjunction with the simulation results presented later show that strategies that fit within the framework of Problem~\ref{prob:budget} exist that reduce $ \avginf $ on average.

An important assumption used herein is that the number of red balls to be added $ \Delta_{r,i}(n) $ is known at least one step ahead of time, so that in particular $ \Delta_{r,i}(n) $ is almost surely constant given $ \F_{n-1} $. A sufficient, but not necessary, condition to satisfy this assumption is for $ \{\Delta_{r,i}(n)\}_{n=1}^{\infty} $ to be set, for all $ i \in V $, before the process begins.


\begin{theorem}\longthmtitle{Individual Urn Proportion Categories}\label{thm:U_n_threshold}
In a general network $ \G = (V,\E) $, if we choose $\{\Delta_{b,i}(n)\}_{n=1}^{\infty}$ so that
\[
	\Delta_{b,i}(n) \geq \frac{\Delta_{r,i}(n)(1 - U_{i,n-1})S_{i,n-1}}{U_{i,n-1}(1-S_{i,n-1})}
\]
almost surely for all $n \in \integerspositive $ and $ i \in V $ (resp. equal to, less than, or equal to) then $ \{U_{i,n}\}_{n=1}^{\infty} $ is a supermartingale (resp. martingale, submartingale) with respect to the natural filtration $ \filt $, i.e.,
\[
	E[{U}_{i,n} | \mathcal{F}_{n-1}] \leq {U}_{i,n-1} \qquad \text{almost surely } \forall n \in \integerspositive.
\]
\end{theorem}

\begin{corollary}\longthmtitle{Network Susceptibility Supermartingale}\label{cor:susceptibility_sup}
In a general network $\G = (V,\E)$, if the curing policies $\{\Delta_{b,i}(t)\}_{t=1}^{\infty}$ obey the bound
\[
	\Delta_{b,i}(n) \geq \frac{\Delta_{r,i}(n)(1 - U_{i,n-1})S_{i,n-1}}{U_{i,n-1}(1-S_{i,n-1})}
\]	
almost surely for all nodes $i \in V$, then the network susceptibility $ \{\tilde{U}_{n}\}_{n=1}^{\infty} $, where $\suscept = \frac{1}{N}\sum_{i=1}^N U_{i,n}$, is a supermartingale with respect to the natural filtration $\{\mathcal{F}_n\}_{n=1}^{\infty}$, i.e.,
\[
	E[\tilde{U}_{n} | \mathcal{F}_{n-1}] \leq \tilde{U}_{n-1} \qquad \text{almost surely } \forall n \in \integerspositive.
\]
\end{corollary}

The proof for Theorem~\ref{thm:U_n_threshold} is provided in Appendix~\ref{app:proof_U_n}. While Corollary~\ref{cor:susceptibility_sup} is useful, the network exposure $ \exposure $ is more closely related to the average infection rate $ \avginf $ than the network susceptibility $ \suscept $, since our draws are taken from the super urn. It is with this in mind that we show the next results, which give us sufficient conditions for $ \{S_{i,n}\}_{n=1}^{\infty} $ and $ \{\tilde{S}_{n}\}_{n=1}^{\infty} $ to be supermartingales.

\begin{theorem}\longthmtitle{Super Urn Proportion Categories}\label{thm:sup_mart}
In a general network $\G = (V,\E)$, if the curing policy $\{\Delta_{b,i}(t)\}_{t=1}^{\infty}$ obeys the lower bound
\begin{equation}\label{eq:sup_lower_bound}
	\Delta_{b,i}(n) > \Delta_{r,i}(n)\frac{S_{i,n-1}}{1-S_{i,n-1}}\max_{k \text{ s.t. } i \in \N_k^{'}} \frac{1-S_{k,n-1}}{S_{k,n-1}} \tag{B1}
\end{equation}	
almost surely for all nodes $i \in V$, then the neighbourhood proportions of red balls $ \{S_{i,n}\}_{n=1}^{\infty} $ are strict supermartingales with respect to the natural filtration $\{\mathcal{F}_n\}_{n=1}^{\infty}$, i.e.
\[
	E[S_{i,n} | \mathcal{F}_{n-1}] < S_{i,n-1} \qquad \text{almost surely } \forall i\in V,n \in \integerspositive.
\]
Furthermore, if the curing policy $\{\Delta_{b,i}(t)\}_{t=1}^{\infty}$ obeys the upper bound
\begin{equation}\label{eq:sup_upper_bound}
	\Delta_{b,i}(n) < \Delta_{r,i}(n)\frac{S_{i,n-1}}{1-S_{i,n-1}}\min_{k \text{ s.t. } i \in \N_k^{'}} \frac{1-S_{k,n-1}}{S_{k,n-1}} \tag{B2}
\end{equation}
almost surely for all nodes $i \in V$, then the neighbourhood proportions of red balls $ \{S_{i,n}\}_{n=1}^{\infty} $ are strict submartingales with respect to the natural filtration $\{\mathcal{F}_n\}_{n=1}^{\infty}$.
\end{theorem}

\begin{corollary}\longthmtitle{Network Exposure Categories}\label{cor:exposure_sup}
In a general network $\G = (V,\E)$, if the curing policies $\{\Delta_{b,i}(t)\}_{t=1}^{\infty}$ obey the lower bound~\eqref{eq:sup_lower_bound} almost surely for all nodes $i \in V$, then the network exposure $ \{\tilde{S}_{n}\}_{n=1}^{\infty} $, where $\exposure = \frac{1}{N}\sum_{i=1}^N S_{i,n}$, is a strict supermartingale with respect to the natural filtration $\{\mathcal{F}_n\}_{n=1}^{\infty}$, i.e.,
\[
	E[\tilde{S}_{n} | \mathcal{F}_{n-1}] < \tilde{S}_{n-1} \qquad \text{almost surely } \forall n \in \integerspositive.
\]
Furthermore, if the curing policies $\{\Delta_{b,i}(t)\}_{t=1}^{\infty}$ obey the upper bound~\eqref{eq:sup_upper_bound} almost surely for all nodes $i \in V$, then the network exposure $ \{\tilde{S}_{n}\}_{n=1}^{\infty} $ is a strict submartingale with respect to the natural filtration $\{\mathcal{F}_n\}_{n=1}^{\infty}$.
\end{corollary}
The proof of Theorem~\ref{thm:sup_mart} is presented in Appendix~\ref{app:proof_sup}. While the duality of these bounds is interesting, in the context of curing we will focus on the lower bound~\eqref{eq:sup_lower_bound}. It is important to note that the policy for $\{\Delta_{b,i}(t)\}_{t=1}^{\infty}$ used in Theorem~\ref{thm:sup_mart} is not a tight lower bound on the curing resources which guarantee that the processes $ \{S_{i,n}\}_{n=1}^{\infty} $ are supermartingales, and hence it is possible that less costly policies exist that still guarantee this property. In particular, strategies may exist which obey the fixed budget $ \budget $ on the amount of curing resources that may be used. However, these results motivate the fact that the search for better policies makes sense, since we know that policies exist that will fight the infection and reduce it on average.

\section{Controlling Epidemics: Algorithmic Strategies}
\label{sec:curing}

The supermartingale results established in the previous section demonstrate the feasibility of a relaxed version of Problem~\ref{prob:budget}, with no budget limitation. In this section, we establish numerical methods to find control policies that find efficient sub-optimal policies for Problem~\ref{prob:budget}, under budget constraints and with having computational complexity in mind. We compare these strategies with the ones obtained from our supermartingale results. A summary of all strategies that will be discussed in this section is given in Table~I.


\begin{table}[ht!]\label{tab:strats}
	\centering
		\ra{1.2}	
	\begin{tabular}{cc}\hline\hline			\\
				$ \stratUn $ & Forcing all $ U_{i,n} $ to be supermartingales (Theorem~\ref{thm:U_n_threshold}): \\
			  & $ \Delta_{bi,}(t) = \frac{\Delta_{r,i}(n)(1 - U_{i,n-1})S_{i,n-1}}{U_{i,n-1}(1-S_{i,n-1})} $
	\\
	\\
\hline
	\\
		$ \stratSn $ & Forcing all $ S_{i,n} $ to be supermartingales (Theorem~\ref{thm:sup_mart}): \\
			 & $ \Delta_{b,i}(t) = \Delta_{r,i}(n)\frac{S_{i,n-1}}{1-S_{i,n-1}}\max_{k \text{ s.t. } i \in \N_k^{'}} \frac{1-S_{k,n-1}}{S_{k,n-1}} $
	\\
	\\
	\hline
	\\
		$ \stratgf $ & Constrained gradient descent algorithm on a simplex: \\
			& Find $ \Delta_{b,i}(t) $ using Algorithm~\ref{alg:grad_flow}
	\\
	\\
	\hline
	\\
		$ \stratheur $ & Ratio of degree, closeness centrality and super urn proportion: \\
			& $ \Delta_{b,i}(t) = \budget\frac{|\N_i| C_i S_{i,t-1}}{\sum_{j=1}^N |\N_j|C_j S_{j,t-1}} $
	\\
	\\
	\hline
	\\
		$ \stratunif $ & Uniformly allocate the budget to all nodes in the network: \\
			& $ \Delta_{b,i}(t) = \frac{\budget}{N} $\\
	\\
\hline\hline			\\
	\end{tabular}
	\caption{Curing Strategies}
\end{table}

Before we present these strategies in details, let us describe briefly how we have evaluated their performance. The simulation platform for these strategies is outlined in Algorithm~\ref{alg:sim}. To achieve comparable results, independent trials of the process are ran with the same initial conditions $ \vec{R} = (R_1,\ldots,R_N) $, $ \vec{B} = (B_1,\ldots,B_N) $, and $  \vec{\Delta}_r = (\Delta_{r,1},\ldots,\Delta_{r,N}) $, for each curing strategy. The results for each strategy is then averaged to evaluate the expected performance. The full simulation results, along with discussions of their implications, are presented in Section~\ref{sec:simulations}.

\begin{algorithm}[!ht]
\caption{Network contagion curing}
\label{alg:sim}
  \begin{algorithmic}
  	\State $ A \gets $ adjacency matrix of the network
	\State $ \vec{R}, \vec{B}, \vec{\Delta}_r \sim \ceil{\mathbf{Uniform}((0,10])}^N $
	\State $ \budget \gets $ budget, $ \Nsum \Delta_{r,i} $
	\State $ C \gets $ number of cases, each with a \textit{strategy}
	\State $ T \gets $ number of trials to run for each case
	\State $ K \gets $ number of time steps for each trial
  	\For{c = 1 : C}
		\For{t = 1 : T}
			\State $\vec{Z}_{c,t} \gets $ \textsc{RunTrial}($A, \vec{R}, \vec{B}, \vec{\Delta}_r, \budget, K, \textit{strategy}$)
		\EndFor
		\State $\vec{Z}_{c} = \frac{1}{T}\sum_{t=1}^{T}\vec{Z}_{c,t}$ 
	\EndFor
  	\\
	\Procedure{RunTrial}{$A, \vec{R}, \vec{B}, \vec{\Delta}_r, \budget, K, \textit{strategy}$}
		\State Initialize $ S_{i,0}, U_{i,0} $ using $ R_i $ and $ B_i $ for all $ i \in V $
		\For{s = 1 : K}
			\State Assign $ \Delta_{b,i}(s) $ using \textit{strategy}, obeying $ \budget $ if required 
			\State Generate $ \vec{Y} \sim \mathbf{Uniform}([0,1])^{N} $
			\If{$ Y_i \leq S_{i,s-1} $}
				\State $ Z_{i,s} = 1 $
			\Else
				\State $ Z_{i,s} = 0 $
			\EndIf
			\State Update $ S_{i,s}, U_{i,s} $ using $ \Delta_{r,i} $, $ \Delta_{b,i}(s) $ for all $ i \in V $
			\State based on $ A $
		\EndFor
	\EndProcedure
  \end{algorithmic}
\end{algorithm}

\subsection{Supermartingale Strategies}\label{subsec:supermart_strats}

The supermartingales results that we have obtained in Section~\ref{sec:control_analytical}, specifically  
Theorems~\ref{thm:U_n_threshold} and~\ref{thm:sup_mart}, naturally lead to a class of curing strategies.
In particular, these strategies guarantee that $ \suscept $ and $ \exposure $, respectively, are supermartingales. It is worth reminding that our theoretical results do not necessarily  imply that average infection rate $ \avginf $ forms a supermartingale. In spite of this, these strategies are still viable options for curing, as far as enough resources are available. We next describe the differences between the strategy given by individual urn proportions, and the one given by super urn proportions.

By Corollary~\ref{cor:susceptibility_sup}, we know that strategy $ \stratUn $ guarantees that the network susceptibility $ \suscept $ will be a supermartingale. Hence we set the curing strategy for each node so that it will force its own individual urn proportion of red balls to be a supermartingale. Since draws are taken form the super urns and not the individual urns, the relationship between the reduction of $ \suscept $ and $ \avginf $ is not a strong one and our simulations suggests that this strategy does not appear to offer a large reduction in the average infection rate in general. In contrast, the curing strategy given by Corollary~\ref{cor:exposure_sup}, where we choose our curing strategy to force the super urn proportions of red balls to be supermartingales for all nodes, performs reasonably well.

%

While these strategies guarantee a reduction in their respective measures, they use an arbitrary amount of curing resources to do so in general. In fact, as we will see later, these strategies always use a large amount of curing resources relative to the impact they have on reducing the average infection rate. To stay within the framework of Problem~\ref{prob:budget}, we will now examine a numerical curing strategy that obeys a fixed budget on the per-step curing resources.

\subsection{Gradient Flow Methods}\label{subsec:grad_flow}

In this section, we employ the well-known gradient descent algorithm~\cite{DPB:95} for Problem~\ref{prob:budget}. As discussed earlier, using $ \avginf $ as a measure of infection is computationally expensive, and hence we instead focus on the network exposure $ \exposure $. While our suggested gradient descent algorithm will not provide the exact answer to Problem~\ref{prob:budget} for reducing $ \avginf$, we will show that it is guaranteed to provide the optimal policy to reduce the closely related measure $ \exposure $.

In Problem~\ref{prob:budget}, our curing policy is constrained by a budget $ \budget $ at each time step and so the feasible set, or set of valid curing policies, for our gradient descent is all policies which do not exceed $ \budget $. However, any optimal policy will make use of the whole budget, and so we consider our feasible set to be $ \X = \left\{\{\Delta_{b,i}(n)\}_{i=1}^N \in \realnonnegative^N \ | \ \Nsum \Delta_{b,i}(n) = \budget \right\} $. Proposition~\ref{prop:grad_flow} shows that for arbitrary initial conditions and network topologies, the problem under study for $ E[\exposure | \F_{n-1}] $ is convex. 

\begin{proposition}\longthmtitle{Gradient descent conditions are met}\label{prop:grad_flow}
In a general network $ \G = (V,\E) $ with arbitrary initial conditions, the expected network exposure $ E[\tilde{S}_{n} | \F_{n-1}] $ is convex with respect to the curing parameters $ \{\Delta_{b,i}(n)\}_{i=1}^N $ for all $ n $. Furthermore, the feasible set
\[
	\X = \left\{\{\Delta_{b,i}(n)\}_{i=1}^N \in \realnonnegative^N \ \Bigg| \ \Nsum \Delta_{b,i}(n) = \budget \right\}
\]
is convex and compact.
\end{proposition}

\begin{figure*}[!th]
\centering
	\subfigure[Facebook group network with 1,363 nodes and 2,425 edges.]{ \includegraphics[width=0.45\linewidth]{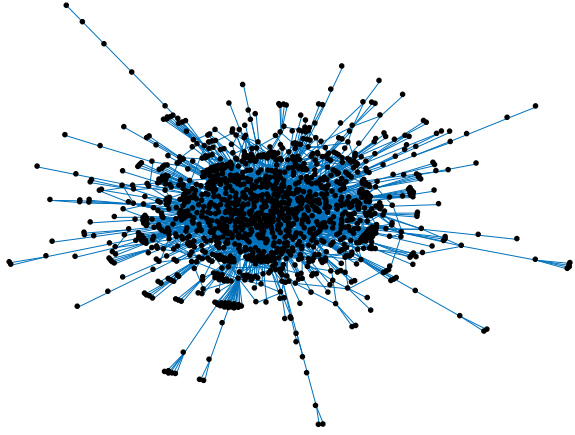} \label{subfig:gf_network}}
\qquad
{
	\psfrag{0}[rr][rr]{\footnotesize{0}}
	\psfrag{0.1}[rr][rr]{\footnotesize{0.1}}
	\psfrag{0.4}[rr][rr]{\footnotesize{0.4}}
	\psfrag{0.5}[rr][rr]{\footnotesize{0.5}}
	\psfrag{250}[rr][rr]{\footnotesize{250}}
	\psfrag{500}[rr][rr]{\footnotesize{500}}
	\psfrag{aa(i)}[rr][rr]{\footnotesize{(i)}}
	\psfrag{aa(ii)}[rr][rr]{\footnotesize{(ii)}}
	\psfrag{aa(iii)}[rr][rr]{\footnotesize{(iii)}}
	\psfrag{aa(iv)}[rr][rr]{\footnotesize{(iv)}}
	\psfrag{aa(v)}[rr][rr]{\footnotesize{(v)}}
	\psfrag{aaAverage r}[rr][rr]{\footnotesize{Average $\rho$}}
	\subfigure[Plot of empirical average infection rate $ \avginf $ compared to the initial level of infection $ \rho = \frac{\Nsum R_i}{\Nsum T_i} $ (lower means less infection).]{ \includegraphics[width=0.45\linewidth]{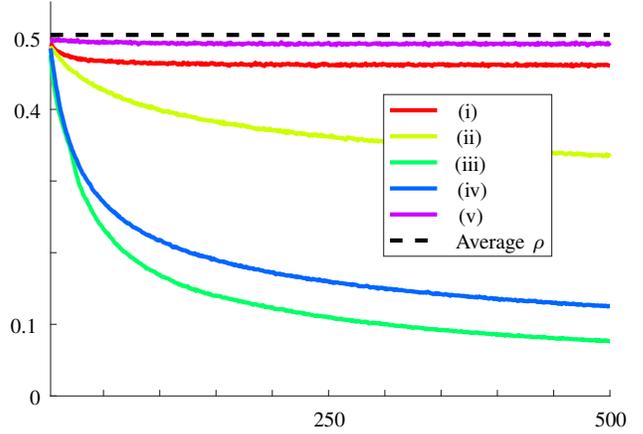} \label{subfig:gf_curing}}
}
\\
{
	\psfrag{250}[rr][rr]{\footnotesize{250}}
	\psfrag{500}[rr][rr]{\footnotesize{500}}
	\psfrag{7500}[rr][rr]{\footnotesize{7,500}}
	\psfrag{8500}[rr][rr]{\footnotesize{8,500}}
	\psfrag{9000}[rr][rr]{\footnotesize{9,000}}		
	\psfrag{11500}[rr][rr]{\footnotesize{11,500}}	
	\psfrag{aa(i)}[rr][rr]{\footnotesize{(i)}}
	\psfrag{aa(ii)}[rr][rr]{\footnotesize{(ii)}}
	\psfrag{aa(iii)}[rr][rr]{\footnotesize{(iii)}}
	\psfrag{aa(iv)}[rr][rr]{\footnotesize{(iv)}}
	\psfrag{aa(v)}[rr][rr]{\footnotesize{(v)}}
	\subfigure[Plot of empirical average usage of curing resources, $ \Nsum \Delta_{b,i}(t) $ (lower means less curing resources used). Note that strategies $ \stratgf $, $ \stratheur $ and $ \stratunif $'s usages are fixed at the per-step budget  $ \budget $, and are overlaid.]{ \includegraphics[width=0.45\linewidth]{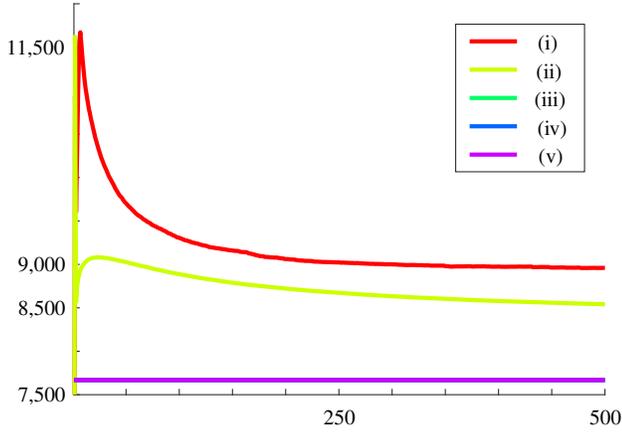} \label{subfig:gf_curing_usage}}
}
\qquad
{
	\psfrag{1000}[rr][rr]{\footnotesize{1,000}}
	\psfrag{3000}[rr][rr]{\footnotesize{3,000}}
	\psfrag{4000}[rr][rr]{\footnotesize{4,000}}	
	\psfrag{6000}[rr][rr]{\footnotesize{6,000}}	
	\psfrag{8000}[rr][rr]{\footnotesize{8,000}}	
	\psfrag{aa(i)}[rr][rr]{\footnotesize{(i)}}
	\psfrag{aa(ii)}[rr][rr]{\footnotesize{(ii)}}
	\psfrag{aa(iii)}[rr][rr]{\footnotesize{(iii)}}
	\psfrag{aa(iv)}[rr][rr]{\footnotesize{(iv)}}
	\psfrag{aa(v)}[rr][rr]{\footnotesize{(v)}}
	\psfrag{aaBudget}[rr][rr]{\footnotesize{Budget}}
	\subfigure[Plot of empirical average wasted curing, $ \Nsum\tsum \Delta_{b,i}(t)Z_{i,t} $ (lower means less curing resources assigned to nodes that did not use them).]{ \includegraphics[width=0.45\linewidth]{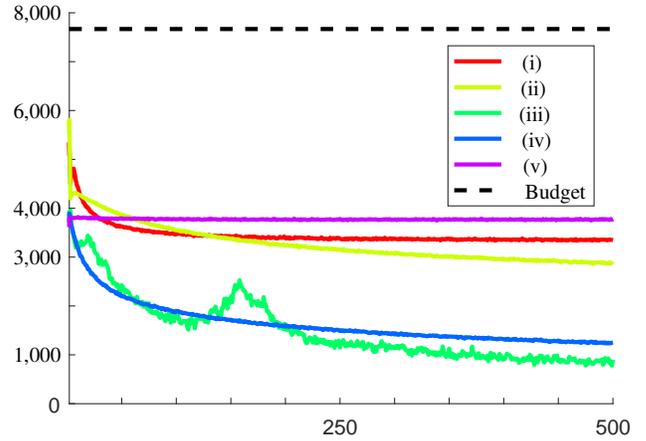} \label{subfig:gf_waste}}
}
\caption{Comparison of all curing strategies presented in Table~I. Simulation results were averaged over 250 trials for each strategy, and altogether took approximately 49 hours on 10 cores of an Intel Xeon processor at 2.20GHz. Initial numbers of balls $R_i$ and $B_i$, and numbers of red balls added $\Delta_{r,i}$ (which remained constant in time), were uniformly randomly assigned for each node but stayed consistent throughout all trials and strategies, while the assignments of $ \{\Delta_{b,i}(t)\}_{t=1}^{\infty} $ were different for each strategy. Since the $\Delta_{r,i}$ are all constant, the budget was set as $ \budget = \Nsum \Delta_{r,i} $.}
\label{fig:gf_sim}
\end{figure*}

The proof for Proposition~\ref{prop:grad_flow} can be found in Appendix~\ref{app:proof_grad}. The structure of the feasible set $ \X $ allows us to employ the simplex constrained gradient descent method, see~\cite[Chapter~2]{DPB:95}; this procedure is fully described in Algorithm~\ref{alg:grad_flow}. The complexity of this algorithm which runs at each time step is $ O(Nsa) $, where $ N $ is the number of nodes, $ s $ is the number of iterations of the gradient descent, and $ \frac{1}{a} $ is the granularity used to find the limit-minimized step size $ \alpha_k $.  While Proposition~\ref{prop:grad_flow} guarantees that the curing policy that this algorithm finds will be optimal for each individual step, it does not guarantee optimality over the entire time horizon. 
In spite of this, as the simulation results in Figure~\ref{fig:gf_sim} show, this curing strategy still outperforms all other curing strategies studied in this paper. The downside of the gradient method is that it is computationally expensive to execute, as it requires intimate knowledge of the state of all nodes in the network. This motivates us to seek other methods which are computationally easier to execute, although they do not perform as well as the gradient descent strategy.

\begin{algorithm}[!ht]
\caption{Constrained gradient descent on a simplex~\cite{DPB:95}}
\label{alg:grad_flow}
  \begin{algorithmic}
	\State $ y_1 = (\budget,0,\ldots,0) $
	\State $ T \gets $ number of iterations to perform
	\\
	\For{$ k = 1 : T $}
		\State $ i = \argmin_{j \in V} \frac{\partial f}{\partial x_j} $
		\State $ [\bar{y}_k]_i = \budget $, and $ [\bar{y}_k]_j = 0 $ for all $ j \neq i $
		\State $ \alpha_k = \argmin_{\alpha \in [0,1]} f(y_k + \alpha(\bar{y}_k - y_k)) $
		\State $ y_{k+1} = y_k +  \alpha_k(\bar{y}_k - y_k) $
	\EndFor
  \end{algorithmic}
\end{algorithm}

\subsection{Heuristic Strategies}\label{subsec:heuristic}

Both sets of strategies identified above come with challenges. The supermartingale strategies are accompanied by analytical results that guarantee that they will improve in general, but they do not obey a fixed budget, nor do they create a significant reduction in the average infection rate. The gradient flow method uses a fixed budget and is provably optimal to reduce the expected network exposure $ E[\exposure | \F_{n-1}] $, but it is computationally costly and requires a large amount of information about the state of infection at every node, including the entire history of draws and values of the curing parameters. As a compromise between these strategies, we present the centrality-infection ratio strategy, which is a heuristic centrality-based strategy designed to allocate the fixed per-step budget $ \budget $.

The idea is to create a ratio to split the budget between all nodes in the network, whose time complexity will be of the order $ O(1) $. We consider three factors when determining how much curing a node should receive: local impact, topological position, and level of infection. Nodes with higher local impact have more neighbours, and hence any healing they receive is immediately distributed to a larger number of nodes. Those with a better topological position are more centrally located within the network, in the sense that the distance from them to all other nodes is smaller. Lastly, nodes with a higher level of infection will need more curing resources to become healthy.

The resulting curing strategy, which we called the centrality-infection ratio, is 
\[
	\Delta_{b,i}(t) = \budget\frac{|\N_i| C_i S_{i,t-1}}{\sum_{j=1}^N |\N_j|C_j S_{j,t-1}}.
\]
To measure local impact of node $ i $, we use the degree, $ | \N_i | $, which measures the number of neighbours for node $ i $. Topological position is determined by calculating the closeness centrality~\cite{AB:50}, which, for node $ i $, is defined as
\[
	C_i := \frac{1}{\sum_{j\in V} d(i,j)},
\]
where $ d(i,j) $ is the length of the shortest path from node $ i $ to node $ j $. Thus $ C_i $ will be higher than $ C_j $ if node $ i $ is closer to all other nodes than node $ j $, in the sense that the paths from node $ i $ will be shorter in total. Hence $ C_i $ gives more importance to the nodes which are more central, and thus have more influence on the overall network. Finally, to measure the level of infection, we use the super urn proportion of red balls $ S_{i,n} $. From~\eqref{eq:cond_eq}, we know that this quantity captures how likely it is for node $ i $ to be infected at this time given the history of the process. Thus we give more importance to nodes who are more likely to be infected, so that we may make them less likely to be infected in the future. This allocation ratio is a generalization of  the best heuristic strategy presented in~\cite{MH-FA-BG:17}, which only used the degree to measure centrality.

The advantage of this heuristic strategy is twofold. Not only does it reduce computational time complexity from $ O(Nsa) $ to $ O(1) $, it is also somewhat distributed in the sense that it does not require constant information from the entire network. Unlike the gradient descent algorithm, strategy $ \stratheur $ simply needs to know information about the network topology and the state of infection of each node. Since we assume that our network's graph is constant in time, this topological information is only required initially and can be used thereafter. The only other information required from the network at large is the sum of the super urn ratios $ \Nsum S_{i,n} $, and hence much less information needs to be communicated through the network for the implementation of this strategy.

Lastly, for comparison reasons we present the uniform curing strategy $ \stratunif $, which splits the budget $ \budget $ equally to all nodes in the network. This provides a benchmark to measure the improvement achieved by more intelligent strategies.

\section{Simulation Results and Discussion}
\label{sec:simulations}

\begin{figure}[!ht]
\centering
{
	\psfrag{0}[rr][rr]{\footnotesize{0}}
	\psfrag{0.1}[rr][rr]{\footnotesize{0.1}}
	\psfrag{0.3}[rr][rr]{\footnotesize{0.3}}
	\psfrag{0.5}[rr][rr]{\footnotesize{0.5}}
	\psfrag{2500}[rr][rr]{\footnotesize{2,500}}
	\psfrag{5000}[rr][rr]{\footnotesize{5,000}}
	\psfrag{aa(i)}[rr][rr]{\footnotesize{(i)}}
	\psfrag{aa(ii)}[rr][rr]{\footnotesize{(ii)}}
	\psfrag{aa(iv)}[rr][rr]{\footnotesize{(iv)}}
	\psfrag{aa(v)}[rr][rr]{\footnotesize{(v)}}
	\psfrag{aaAverage r}[rr][rr]{\footnotesize{Average $\rho$}}
	\includegraphics[width=\linewidth]{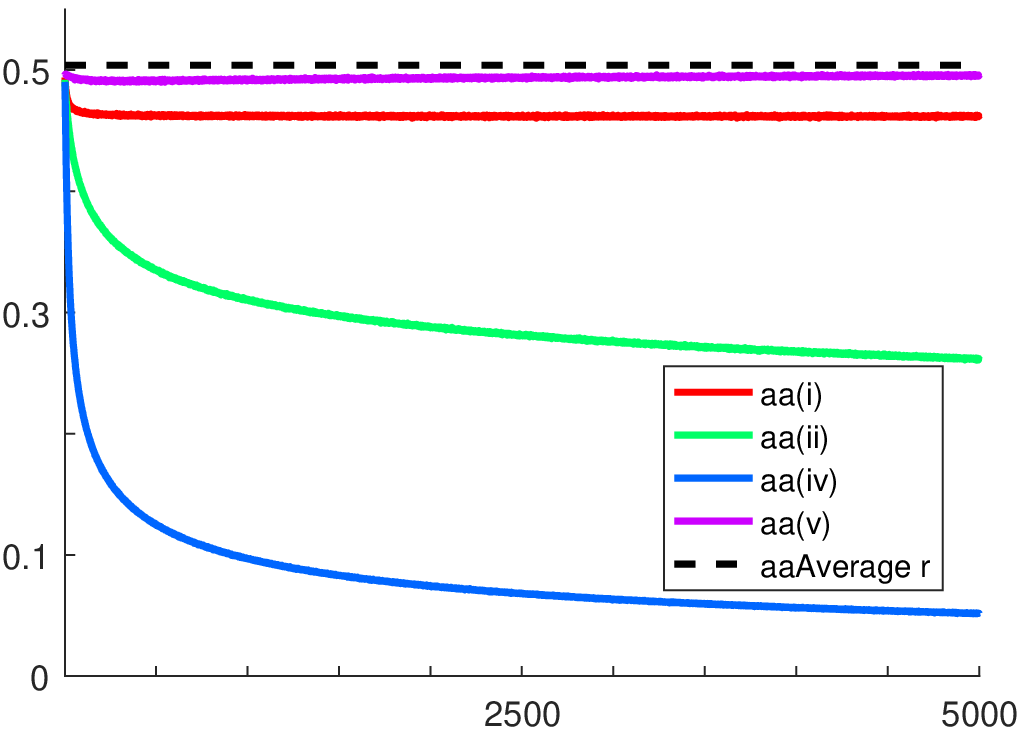} \label{subfig:fb_curing}
}
\caption{Plot of empirical average infection rate $ \avginf $ on the network shown in~\ref{subfig:gf_network} for a longer time frame. Strategies used are shown in Table~I. The simulations presented here were performed identically to those described in Figure~\ref{fig:gf_sim}, with all initial conditions consistent between trials and strategies. Strategies were averaged over 1,000 trials, and altogether took approximately 30 hours on 12 cores of an Intel Xeon processor at 2.20GHz.}
\label{fig:fb_sim}
\end{figure}

\begin{figure*}[th]
\centering
{
	\psfrag{0}[rr][rr]{\footnotesize{0}}
	\psfrag{0.1}[rr][rr]{\footnotesize{0.1}}
	\psfrag{0.4}[rr][rr]{\footnotesize{0.4}}
	\psfrag{0.5}[rr][rr]{\footnotesize{0.5}}
	\psfrag{500}[rr][rr]{\footnotesize{500}}
	\psfrag{1000}[rr][rr]{\footnotesize{1,000}}
	\psfrag{aa(iv)}[rr][rr]{\footnotesize{(iv)}}
	\psfrag{aa(v)}[rr][rr]{\footnotesize{(v)}}
	\psfrag{aaAverage r}[rr][rr]{\footnotesize{Average $ \rho $}}
	\subfigure[Plot of empirical average infection rate $ \avginf $ over $ 1,000 $ time steps.]{ \includegraphics[width=0.45\linewidth]{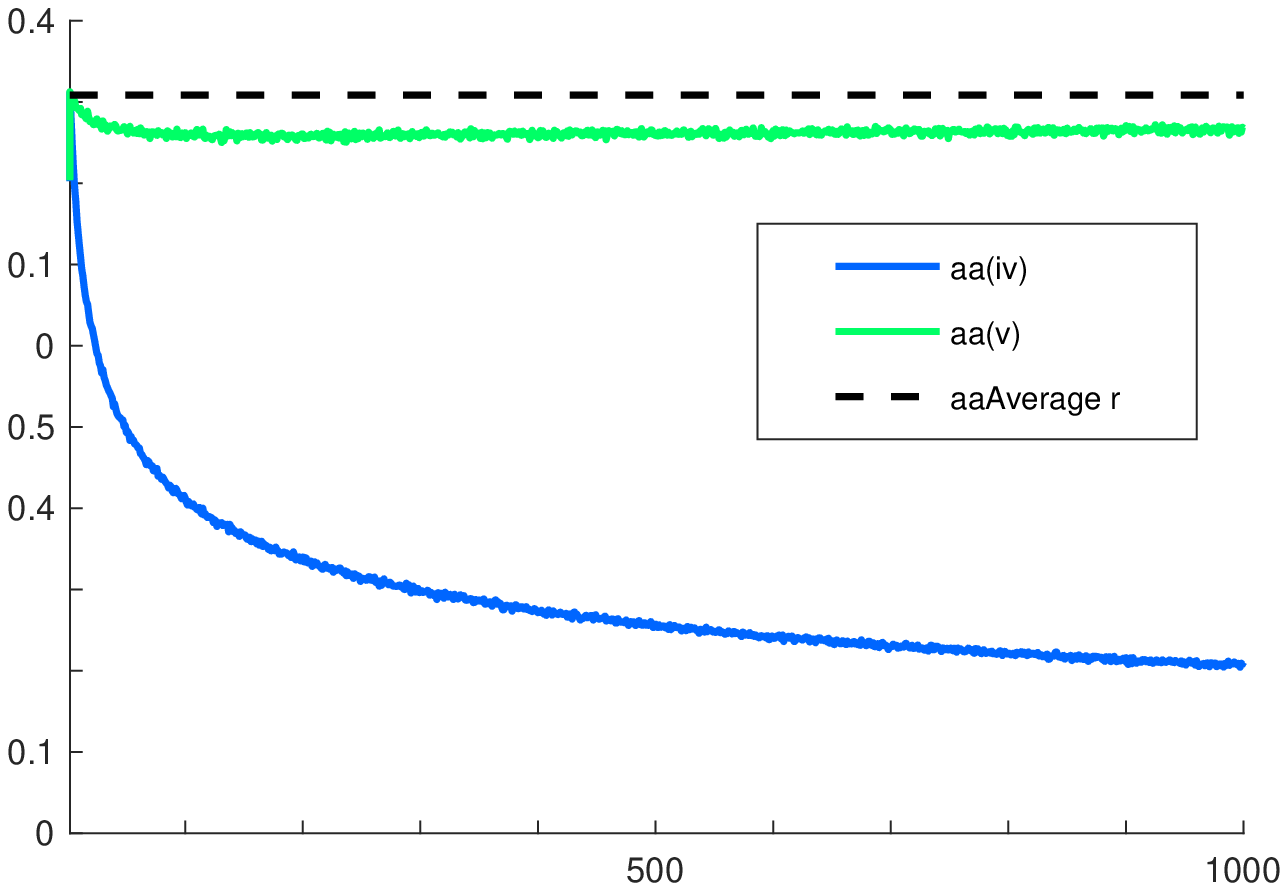} \label{subfig:st_curing}}
}
\qquad
{
	\psfrag{0}[rr][rr]{\footnotesize{0}}
	\psfrag{0.1}[rr][rr]{}
	\psfrag{0.2}[rr][rr]{}
	\psfrag{0.3}[rr][rr]{}
	\psfrag{0.4}[rr][rr]{}
	\psfrag{0.5}[rr][rr]{\footnotesize{0.5}}
	\psfrag{0.6}[rr][rr]{}
	\psfrag{0.7}[rr][rr]{}
	\psfrag{0.8}[rr][rr]{}
	\psfrag{0.9}[rr][rr]{}
	\psfrag{1}[rr][rr]{\footnotesize{1}}
	\subfigure[Plot of initial level of individual infection for each node $ U_{i,0} $.]{ \includegraphics[width=0.45\linewidth]{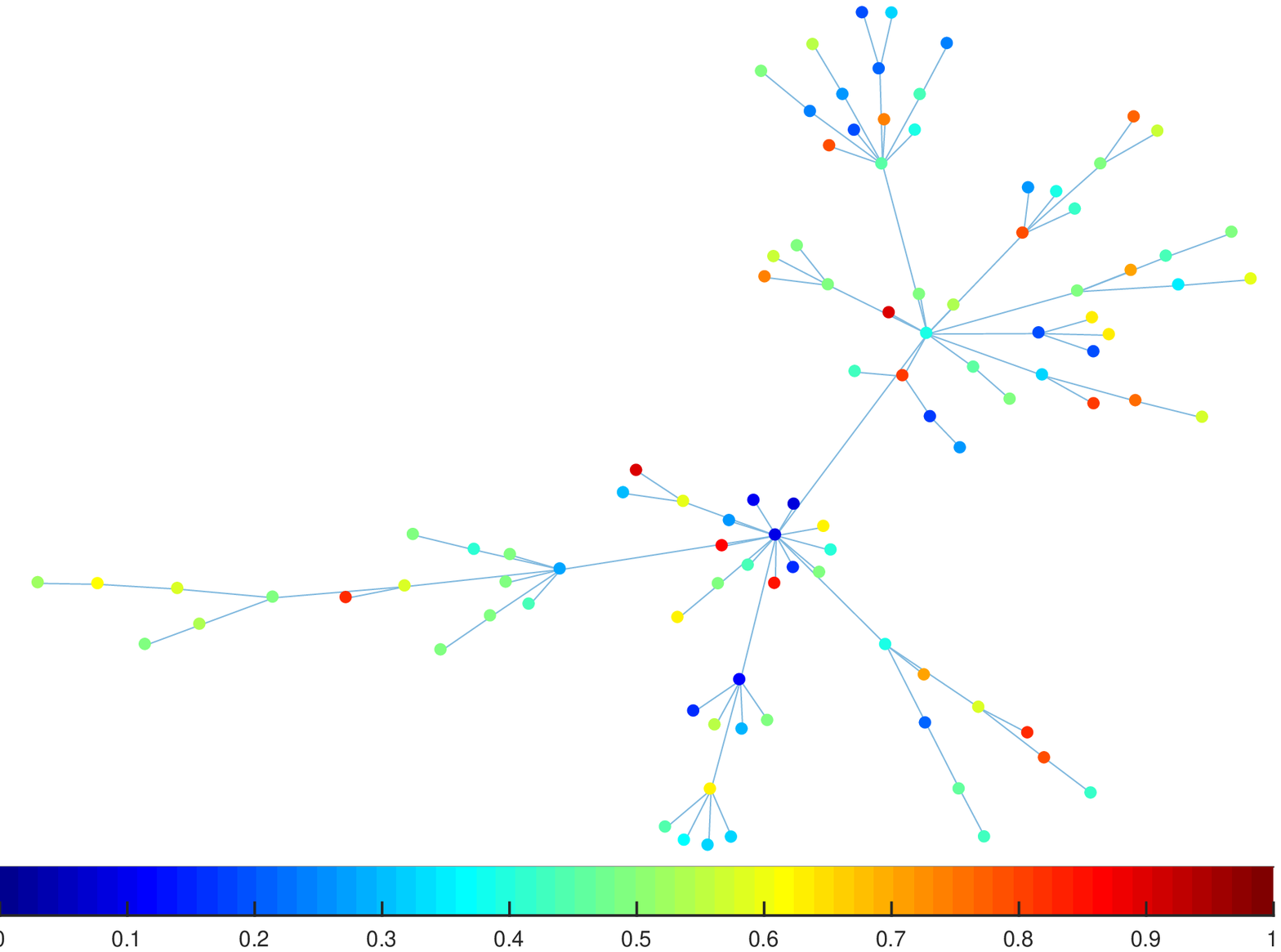} \label{subfig:st_init}}
}
{
	\psfrag{0}[rr][rr]{\footnotesize{0}}
	\psfrag{0.1}[rr][rr]{}
	\psfrag{0.2}[rr][rr]{}
	\psfrag{0.3}[rr][rr]{}
	\psfrag{0.4}[rr][rr]{}
	\psfrag{0.5}[rr][rr]{\footnotesize{0.5}}
	\psfrag{0.6}[rr][rr]{}
	\psfrag{0.7}[rr][rr]{}
	\psfrag{0.8}[rr][rr]{}
	\psfrag{0.9}[rr][rr]{}
	\psfrag{1}[rr][rr]{\footnotesize{1}}
	\subfigure[Plot of final level of individual infection for each node $ U_{i,n} $ using strategy $ \stratheur $, the centrality-infection ratio.]{ \includegraphics[width=0.45\linewidth]{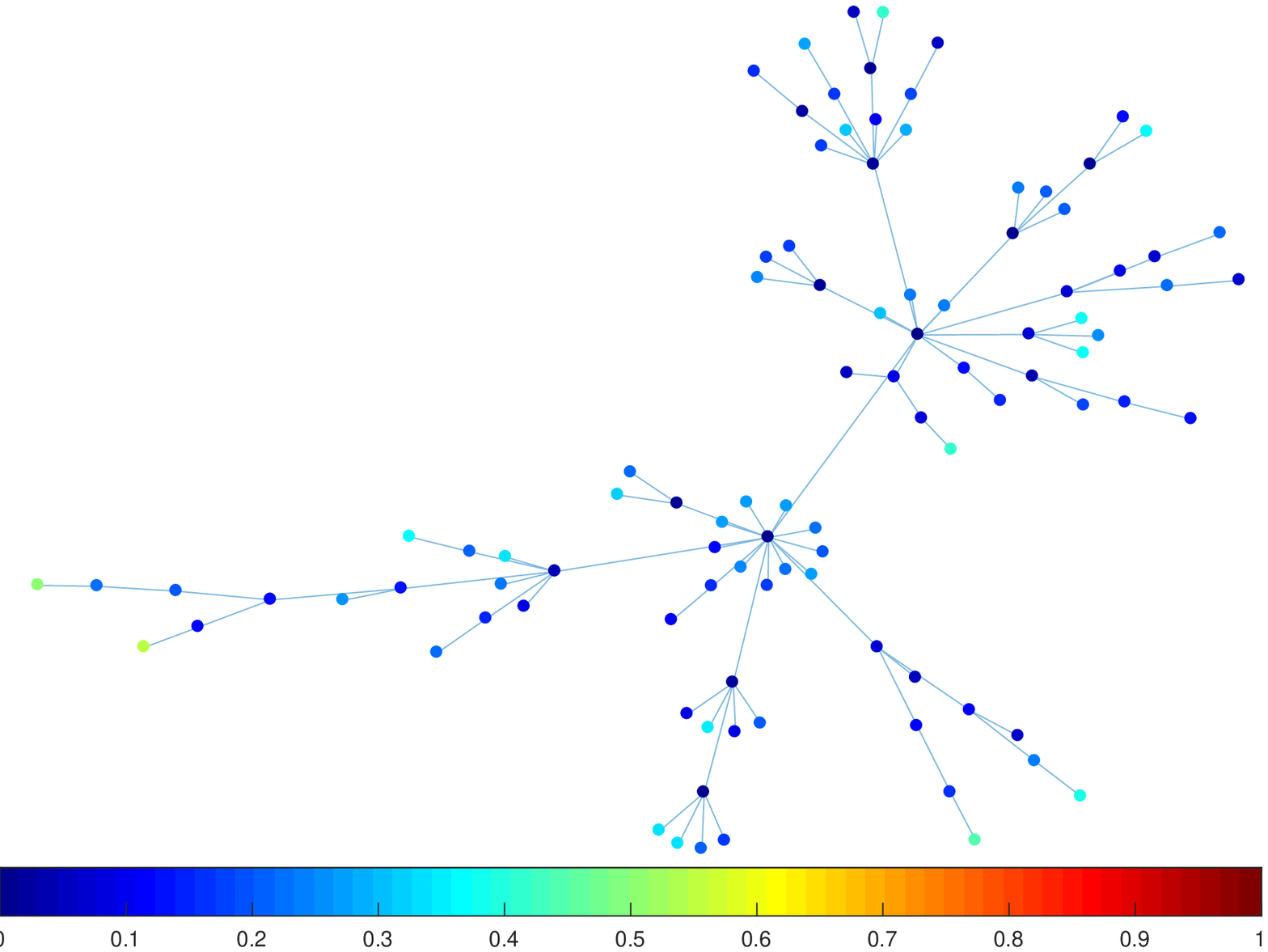} \label{subfig:st_heur}}
}
\quad
{
	\psfrag{0}[rr][rr]{\footnotesize{0}}
	\psfrag{0.1}[rr][rr]{}
	\psfrag{0.2}[rr][rr]{}
	\psfrag{0.3}[rr][rr]{}
	\psfrag{0.4}[rr][rr]{}
	\psfrag{0.5}[rr][rr]{\footnotesize{0.5}}
	\psfrag{0.6}[rr][rr]{}
	\psfrag{0.7}[rr][rr]{}
	\psfrag{0.8}[rr][rr]{}
	\psfrag{0.9}[rr][rr]{}
	\psfrag{1}[rr][rr]{\footnotesize{1}}
	\subfigure[Plot of final level of individual infection for each node $ U_{i,n} $ using strategy $ \stratunif $, uniform curing.]{ \includegraphics[width=0.45\linewidth]{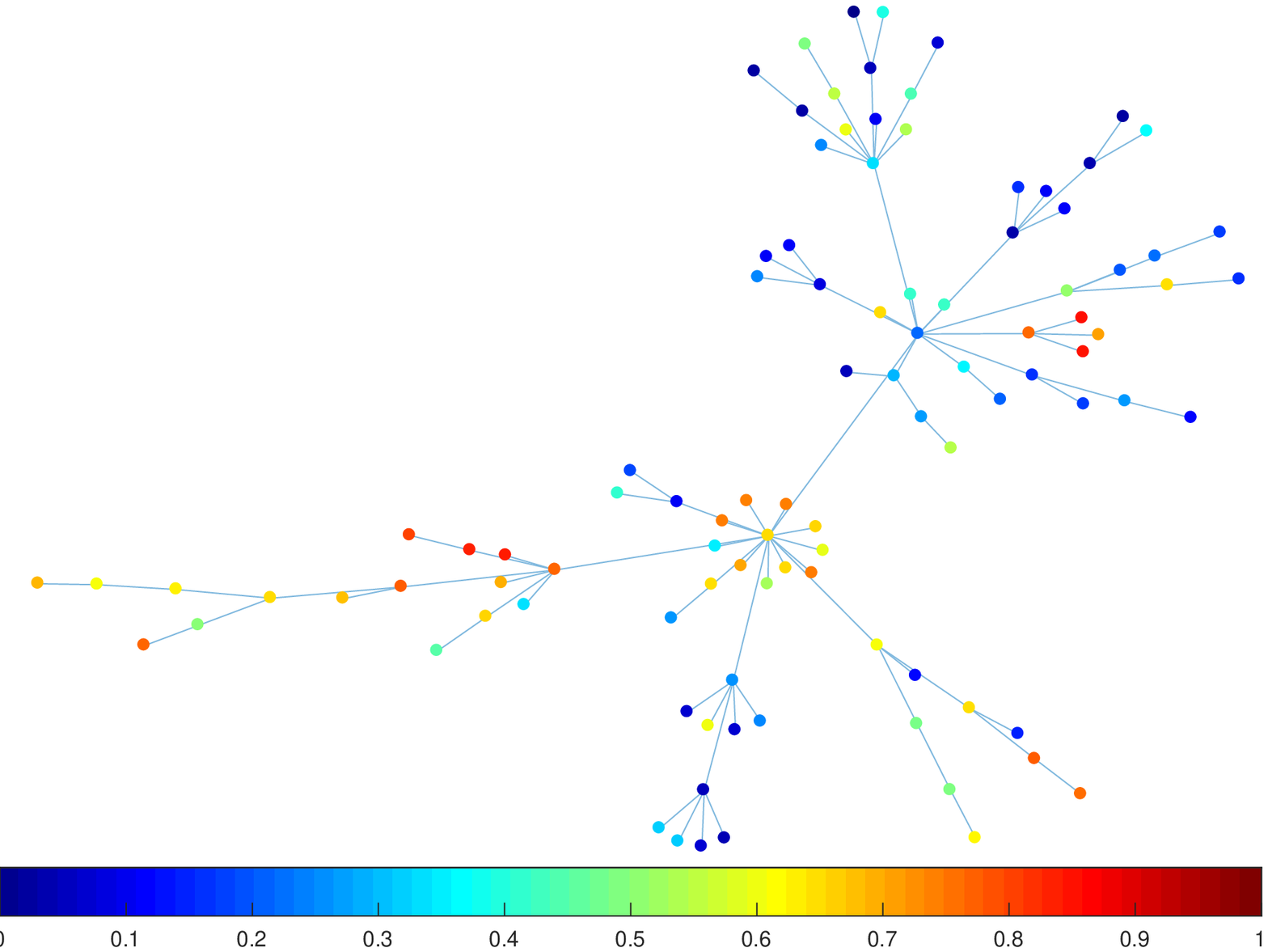} \label{subfig:st_unif}}
}
\caption{Comparison of curing strategies $ \stratheur $ and $ \stratunif $ on a Barabasi-Alberth network~\cite{RA-ALB:02} with 100 nodes and 99 edges. Here blue represents total healthiness ($ U_{i,n} = 0 $) while red represents total infection ($ U_{i,n} = 1 $). Simulation results were averaged over 1000 trials for each strategy, and altogether took approximately 5 minutes on a 4-core Intel Core i7 processor at 2.20 GHz. This simulation was performed identically to those described in Figure~\ref{fig:gf_sim}, with all initial conditions consistent between trials and strategies.}
\label{fig:st_sim}
\end{figure*}

In order to confirm the results of Theorems~\ref{thm:U_n_threshold} and~\ref{thm:sup_mart}, a number of simulations were performed; the pseudocode is outlined in Algorithm~\ref{alg:sim}. While the simulations performed had the numbers of red balls added $ \Delta_{r,i} $ vary between nodes, they were constant in time. This was done to simplify the choice of the per-step budget, and does not affect the execution of the simulations themselves. All initial conditions used in the simulations herein, as well as videos displaying the average performance of the curing strategies, are available online.\footnote{\href{http://bit.ly/2szl8PY}{See: \color{blue}{http://bit.ly/2szl8PY}}}

\subsection{Simulation setup}

The network shown in Figure~\ref{subfig:gf_network} was generated by using a tool~\cite{BR:13} to crawl through 500 posts in a Facebook group. Individuals who created posts or interacted with others' content are represented by nodes, while edges are created if individuals interacted with the post or comment of another (by commenting on the post, or liking the post or comment). The resulting graph has 1,363 nodes and 2,425 edges, and by design represents the topology of a real social network.

We now provide a detailed description of the simulation, as described in Algorithm~\ref{alg:sim}. The values of $ R_i $, $ B_i $ and $ \Delta_{r,i} $ were uniformly randomly assigned for each node as integers between 1 and 10. These values remained consistent for all strategies and throughout all trials that were performed. Since the values for $ \Delta_{r,i} $ were fixed over time, the per-step budget was set at $ \budget = \Nsum \Delta_{r,i} $. With the initial conditions set, a number of trials were performed for each strategy. Each trial was performed by successively drawing balls from super urns for a fixed number of time steps. At time $ t $, we first assigned the curing $ \Delta_{b,i}(t) $ based on the strategy selected. Then a uniform random variable on $ [0,1] $, $ Y_i $, was generated for each node $ i $ and compared to the super urn proportion. If $ Y_i < S_{i,t-1} $ then we say that a red ball was drawn and so $ Z_{i,t} = 1 $, otherwise we drew black and so $ Z_{i,t} = 0 $. Based on what was drawn, we added $ \Delta_{r,i} $ red or $ \Delta_{b,i}(t) $ black balls into node $ i $'s urn, and hence its super urn and those of its neighbours. At the end of each trial the draw variables were saved, and then averaged over all trials to produce the empirical performance of the curing strategy.

\subsection{Discussion of Simulation Results}

The comparisons of all strategies outlined in Section~\ref{sec:curing} can be seen in Figures~\ref{fig:gf_sim} and~\ref{fig:fb_sim}. It is important to note that only strategies $ \stratgf $, $ \stratheur $ and $ \stratunif $ in Table~I have a budget $ \budget $ on the amount of curing they can use, and the other two strategies vary the total curing they use in time; the amount of resources each strategy consumes is shown in Figure~\ref{subfig:gf_curing_usage}. Figure~\ref{subfig:gf_waste} displays the average wasted curing resources for each strategy, which we will define later.

Figures~\ref{subfig:gf_curing} and~\ref{fig:fb_sim} compare the performance of all strategies described in Section~\ref{sec:curing} on a Facebook network. Figure~\ref{subfig:gf_curing} includes the gradient flow algorithm, while Figure~\ref{fig:fb_sim} shows all other strategies over a longer time horizon. The benchmark uniform strategy $ \stratunif $ performs the worst, which is to be expected. Although $ \stratgf $ is only proven to be optimal for the expected network exposure $ E[\exposure | \F_{n-1}] $, we observe that it is effective for the average infection rate $ \avginf $ as well; as previously mentioned, this strategy outperforms all other curing strategies described in this paper. However, the heuristic strategy $ \stratheur $ shows similar performance with dramatic improvements in computational complexity. The supermartingale strategies $ \stratUn $ and $ \stratSn $ both reduce $ \avginf $ below the initial average infection rate in the network $ \rho $, but are less effective in doing so than the other two methods. Strategy $ \stratUn $ sees only an immediate small reduction in $ \avginf $, while strategy $ \stratSn $ continuously decreases $ \avginf $. This empirical evidence follows~\eqref{rmk:measures}, and further supports our assertion that we should focus on the network exposure $ \exposure $ to reduce the average infection rate.

In Figure~\ref{subfig:gf_curing_usage} we examine the amount of curing resources used by each strategy. Since strategies $ \stratgf $, $\stratheur $ and $ \stratunif $ all obey a per-step budget constraint their usages are fixed. Both supermartingale strategies, which may use arbitrary amounts of curing resources, initially use a larger amount of curing resources and then reduce their usage. This initial expenditure is the cost required to turn the measures into supermartingales, after which the strategies only need to maintain the property which requires less resources. Strategy $ \stratUn $'s usage appears to decay exponentially to an almost constant amount, while strategy $ \stratSn $ continues to decrease its usage in time. Further, strategy $ \stratUn $ uses almost $ 50\% $ more curing resources than the budget $ \budget $ initially, while strategy $ \stratSn $'s initial usage is only around $ 18\% $ higher than $ \budget $. This is likely because $ \stratUn $ is a selfish strategy; it only considers what is happening on an individual node level. In contrast, strategy $ \stratSn $ considers the infection in local neighbourhoods of nodes, and hence it is more judicious with its application of resources to specific nodes.

The amount of curing resources wasted by each strategy is displayed in Figure~\ref{subfig:gf_waste}. Waste is defined as curing resources which were assigned to nodes that did not use them since they displayed ``infected'' behaviour at that time, and is hence computed as $ \Nsum\tsum \Delta_{b,i}(t)Z_{i,t} $. We observe an intuitive correlation between the amount of resources wasted and curing performance: strategies which waste less resources tend to be more effective at reducing the average infection rate $ \avginf $. However, this does not tell the full story. The gradient flow algorithm has several spikes where it wastes more resources than the centrality-infection ratio $ \stratheur $, but this does not appear to affect its curing performance. These likely occur because the gradient descent tends to focus on clusters of points when assigning curing, and hence the rest of the network becomes more imfected. Then, when it switches to another cluster, the nodes are more infected and hence it wastes more resources until it sufficiently cures the infection in that cluster. Strategy $ \stratUn $ initially wastes less than strategy $ \stratSn $ even though it uses more curing resources, and it still performs worse with respect to reduction in $ \avginf $. This suggests that optimal curing strategies do not simply waste less, but also intelligently allocate their curing resources to make the best use of them.

Figure~\ref{fig:st_sim} shows the initial and final state of infection of all nodes in a randomly generated network for two different curing strategies. We define the state of infection for node $ i $ at time $ t $ by its individual proportion of red balls $ U_{i,t} $. Here we use a small Barabasi-Albert network so the states at the node level are more visible. Barabasi-Albert networks are randomly generated through preferential attachment and are widely used in the literature, as they have been shown to exhibit the properties of real social networks~\cite{RA-ALB:02}. In Figure~\ref{subfig:st_curing} we see that for such a network the centrality-infection ratio $ \stratheur $ dramatically outperforms uniform curing $ \stratunif $, as was the case for the social network shown in Figure~\ref{subfig:gf_network}. After 1,000 time steps, strategy $ \stratheur $ reduced the average infection rate $ \avginf $ to about $ 15\% $, and no node is worse off than it was to begin with. In contrast, strategy $ \stratunif $ barely reduced $ \avginf $ below the initial average infection rate $ \rho $, and the individual infection of some nodes rose above $ 90\% $. This result illustrates the fact that intelligent allocation of curing resources is not only important to reduce the network-wide average infection rate, but the infection of individual nodes as well.

\section{Conclusion and Future Work}\label{sec:conclusion}

In this paper we examined the problem of curing epidemics using a network contagion model adapted from the Polya urn process. We formulated an optimal control problem and provided analytical results that showed that finding solutions is a worthwhile endeavour. We used theoretical, numerical and heuristic curing strategies to attempt to cure the epidemic, and evaluated their performance using simulations.

Future work with this model could include the statement of different curing problems. A budget could be assigned over a finite time horizon instead of on a per-step basis, and strategies would need to judiciously use this limited supply to reduce infection as much as possible in the time window. The problem could even be reversed, so that some desirable healthiness threshold is provided and one could find the lowest possible budget that would guarantee that the average infection rate would be at or below the threshold. Such a problem could be examined for a per-step or fixed horizon budget.

\begin{appendix}

\subsection{Proof of Theorem~\ref{thm:U_n_threshold}}\label{app:proof_U_n}
This result is a generalization of Theorem 4.6 in~\cite{MH-FA-BG:17-2}, since here we allow $ \{\Delta_{r,i}(t)\}_{t=1}^{\infty} $ to vary in time. As such, some minor steps are omitted.

We will start with the case of a supermartingale. That is, we wish to show that almost surely for all $n \in \integerspositive $,
\[
	E[U_{i,n} \ | \ \F_{n-1}] - U_{i,n-1} \leq 0 \Leftrightarrow E[U_{i,n} - U_{i,n-1} \ | \ \F_{n-1}] \leq 0,
\]
since $ U_{i,n-1} $ is almost surely constant given $ \F_{n-1} $. Take $ X_{i,n} $ as in~\eqref{eq:X_n}. We then compute the difference $ U_{i,n} - U_{i,n-1} $,
\begin{align*}
	&U_{i,n} - U_{i,n-1} \\
	&= \frac{R_i + \sum_{t=1}^n \Delta_{r,i}(t)Z_{i,t}}{X_{i,n}} - \frac{R_i + \sum_{t=1}^{n-1} \Delta_{r,i}(t)Z_{i,t}}{X_{i,n-1}} \\
	&= \frac{\Delta_{r,i}(n) Z_{i,n}}{X_{i,n}} - \frac{(R_i + \sum_{t=1}^{n-1} \Delta_{r,i}(t)Z_{i,t})(X_{i,n} - X_{i,n-1})}{X_{i,n-1}X_{i,n}} \\
	&= \frac{\Delta_{r,i}(n) Z_{i,n} - U_{i,n-1}(\Delta_{r,i}(n) Z_{i,n} +\Delta_{b,i}(n)(1-Z_{i,n}))}{X_{i,n}}.
\end{align*}
Since $ X_{i,n} > 0 $ almost surely, for all $ n \in \integerspositive $, it will not change the sign of the inequality later on, and so we can ignore it to focus only on the numerator. Thus we wish to check if, almost surely,
\begin{align*}
E\big[\Delta_{r,i}(n) Z_{i,n} - U_{i,n-1}&(\Delta_{r,i}(n)Z_{i,n} \\
&\quad + \Delta_{b,i}(n)(1-Z_{i,n})) | \F_{n-1}\big] \leq 0.
\end{align*}
Now if the curing policy $\{\Delta_{b,i}(n)\}_{n=1}^{\infty}$ for node $ i $ satisfies the bound given:
\[
\Delta_{b,i}(n) \geq \frac{\Delta_{r,i}(n)(1 - U_{i,n-1})S_{i,n-1}}{U_{i,n-1}(1-S_{i,n-1})},
\]
then almost surely,
\begin{align*}
&E\big[\Delta_{r,i}(n) Z_{i,n}(1 - U_{i,n-1}) - U_{i,n-1}(1-Z_{i,n})\Delta_{b,i}(n) | \F_{n-1}\big] \\
	&\leq E\Bigg[\Delta_{r,i}(n) Z_{i,n}(1- U_{i,n-1}) - U_{i,n-1}(1-Z_{i,n})\\
	&\qquad\qquad\times\frac{\Delta_{r,i}(n)(1 - U_{i,n-1})S_{i,n-1}}{U_{i,n-1}(1-S_{i,n-1})} \Bigg| \F_{n-1} \Bigg] \\
	&= \Delta_{r,i}(n)(1- U_{i,n-1})\left[S_{i,n-1} - (1-S_{i,n-1})\frac{S_{i,n-1}}{1-S_{i,n-1}} \right] \\
	&= 0,
\end{align*}
where the second to last equality comes from the fact that $ E[Z_{i,n} | \F_{n-1}] = P(Z_{i,n} = 1 | \F_{n-1}) = S_{i,n-1} $ almost surely by~\eqref{eq:cond_eq}, and that $S_{i,n-1}$ is almost surely constant given $\F_{n-1}$. Thus as long as $ \Delta_{b,i}(n) $ obeys this bound almost surely for all $ n \in \integerspositive $, $ \{U_{i,n}\}_{n=1}^{\infty} $ is a supermartingale with respect to $ \{Z_n\}_{n=1}^{\infty} $. Similarly, if $ \Delta_{b,i}(n) $ is almost surely equal (resp. less than or equal) to this bound, $ \{U_{i,n}\}_{n=1}^{\infty} $ is a martingale (resp. submartingale) with respect to $ \{\F_n\}_{n=1}^{\infty} $. \oprocend

\subsection{Proof of Theorem~\ref{thm:sup_mart}}\label{app:proof_sup}

We will focus on the case of a supermartingale, since the proof for submartingales follows similarly. First, note that the question of $\{S_{i,n}\}_{n=1}^{\infty}$ being a strict supermartingale is equivalent to
\[
	E[S_{i,n} | \mathcal{F}_{n-1}] - S_{i,n-1} < 0
\]
where $\{\mathcal{F}_{n}\}$ is the natural filtration (indeed, we can just condition on $Z^{n-1}$). Note, in particular, that $E[Z_{i,t} | \F_{n}] = Z_{i,t}$ almost surely, for all $i \in V$ and $t \in \{1,\ldots,n\}$, and the same is true for $\{S_{i,t}\}_{t=1}^n$. Then almost surely, as in Theorem~\ref{thm:U_n_threshold},
\begin{align*}
	S_{i,n} - S_{i,n-1} 
	&= \frac{S_{i,n-1}(\bar{X}_{i,n-1} - \bar{X}_{i,n}) + \bar{Z}_{r,i,n}}{\bar{X}_{i,n}}.
\end{align*}
Since $ \bar{X}_{i,n} > 0 $ almost surely for all $ n \in \integerspositive $ and all $ i \in V $, we can ignore it. Further, since $ S_{i,n-1} $ is almost surely constant, we need to only check if
\begin{align*}
	&E\left[S_{i,n-1}(\bar{X}_{i,n-1} - \bar{X}_{i,n}) + \bar{Z}_{r,i,n} \ | \ \F_{n-1}\right] < 0 \\
	\Leftrightarrow &E\left[(1 - S_{i,n-1})\bar{Z}_{r,i,n} - S_{i,n-1}\bar{Z}_{b,i,n} \ | \ \F_{n-1}\right] < 0,
\end{align*}
since $ \bar{X}_{i,n-1} - \bar{X}_{i,n} = -\bar{Z}_{r,i,n} - \bar{Z}_{b,i,n} $. Now let the lower bound~\eqref{eq:sup_lower_bound} be satisfied:
\begin{align*}
	\Delta_{b,j}(n) > \Delta_{r,j}(n)\frac{S_{j,n-1}}{1-S_{j,n-1}}\max_{k \text{ s.t. } j \in \N_k^{'}} \frac{1-S_{k,n-1}}{S_{k,n-1}}.
\end{align*}
Notice $ E[Z_{j,n} | \F_{n-1}] = S_{j,n-1} $ almost surely, so we have
\begin{align*}
	&E\left[S_{i,n-1}(\bar{X}_{i,n-1} - \bar{X}_{i,n}) + \bar{Z}_{r,i,n} \ | \ \F_{n-1}\right] \\
	&< E\Bigg[(1 - S_{i,n-1})\bar{Z}_{r,i,n} - S_{i,n-1}\jsum \Delta_{r,j}(n)\frac{S_{j,n-1}}{1 - S_{j,n-1}} \\
	&\qquad\qquad \times \max_{k \text{ s.t. } j \in \N_k^{'}} \frac{1-S_{k,n-1}}{S_{k,n-1}}(1-Z_{j,n}) \Bigg| \F_{n-1}\Bigg] \\
	&=\jsum \Delta_{r,j}(n)S_{j,n-1}(1 - S_{i,n-1}) - \Delta_{r,j}(n)\frac{S_{j,n-1}}{1 - S_{j,n-1}} \\
	&\qquad\qquad \times S_{i,n-1}\max_{k \text{ s.t. } j \in \N_k^{'}} \frac{1-S_{k,n-1}}{S_{k,n-1}}(1-S_{j,n-1})  \\
	&=\jsum \Delta_{r,j}(n) S_{j,n-1}\left[1 - S_{i,n-1}\,\mathclap{\Bigg(}\, 1 + \quad\enskip\mathclap{\max_{k \text{ s.t. } j \in \N_k^{'}}}\quad\enskip \frac{1-S_{k,n-1}}{S_{k,n-1}}\Bigg)\right] \\
	&=\jsum\Delta_{r,j}(n) S_{j,n-1}\left[1 - \frac{S_{i,n-1}}{\min_{k \text{ s.t. } j \in \N_k^{'}}S_{k,n-1}}\right].
\end{align*}
Now note that $ j \in \N_i' $ and hence, in particular, $ \min_{k \text{ s.t. } j \in \N_k^{'}}S_{k,n-1} \leq S_{i,n-1} $ almost surely, and all other quantities are non-negative. Thus, with our value of $ \Delta_{b,j}(n) $ for all $ j \in \N_i' $, we have almost surely
\begin{align*}
	&E\left[S_{i,n-1}(\bar{X}_{i,n-1} - \bar{X}_{i,n}) + \bar{Z}_{i,n} \ | \ \F_{n-1}\right] \\
	&< \jsum\Delta_{r,j}(n) S_{j,n-1}\left[1 - \frac{S_{i,n-1}}{\min_{k \text{ s.t. } j \in \N_k^{'}}S_{k,n-1}}\right] \\
	& \leq 0.
\end{align*} 
Thus, for any $ i \in V $, if $ \{\Delta_{b,i}(n)\}_{n=1}^{\infty} $ obeys the lower bound~\eqref{eq:sup_lower_bound} almost surely, the neighbourhood proportion of red balls $ \{S_{i,n}\}_{n=1}^{\infty} $ is a strict supermartingale.

For the case of a strict submartingale, notice  if the upper bound~\eqref{eq:sup_upper_bound} is satisfied:
\begin{align*}
	\Delta_{b,j}(n) < \Delta_{r,j}(n)\frac{S_{j,n-1}}{1-S_{j,n-1}}\min_{k \text{ s.t. } j \in \N_k^{'}} \frac{1-S_{k,n-1}}{S_{k,n-1}},
\end{align*}
then similarly to the case of a supermartingale,
\begin{align*}
&E\left[S_{i,n-1}(\bar{X}_{i,n-1} - \bar{X}_{i,n}) + \bar{Z}_{i,n} \ | \ \F_{n-1}\right] \\
	&> E\Bigg[(1 - S_{i,n-1})\bar{Z}_{r,i,n} - S_{i,n-1}\jsum \Delta_{r,j}(n)\frac{S_{j,n-1}}{1 - S_{j,n-1}} \\
	&\qquad\qquad \times \min_{k \text{ s.t. } j \in \N_k^{'}} \frac{1-S_{k,n-1}}{S_{k,n-1}}(1-Z_{j,n}) \Bigg| \F_{n-1}\Bigg] \\
	&=\jsum\Delta_{r,j}(n) S_{j,n-1}\left[1 - \frac{S_{i,n-1}}{\max_{k \text{ s.t. } j \in \N_k^{'}}S_{k,n-1}}\right] \\
	&\geq 0,
\end{align*}
since $ \max_{k \text{ s.t. } j \in \N_k^{'}}S_{k,n-1} \geq S_{i,n-1} $ almost surely. Thus if $ \{\Delta_{b,i}(n)\}_{n=1}^{\infty} $ obeys the upper bound~\eqref{eq:sup_upper_bound} almost surely, the neighbourhood proportion of red balls $ \{S_{i,n}\}_{n=1}^{\infty} $ is a strict submartingale.\oprocend

\subsection{Proof of Proposition~\ref{prop:grad_flow}}\label{app:proof_grad}

First note that as a function of the parameters $ x = (\Delta_{b,1}(n),\ldots,\Delta_{b,N}(n))$, $ E[\tilde{S}_{n} | \F_{n-1}] $ is of the form $ f_n(x) = \frac{1}{N}\Nsum \frac{c_i}{d_i + \sigma_{i}(x)} $, where from~\eqref{eq:S_n}, we can see that
\begin{align*}
	c_i &= \bR_i + \Delta_{r,j}(n)E[Z_{j,n} | \F_{n-1}] + \sum_{t=1}^{n-1}\bar{Z}_{r,i,t}, \\
	d_i &= c_i + \bB_i + \sum_{t=1}^{n-1}\bar{Z}_{b,i,t}, \text{ and} \\
	\sigma_{i}(x) &= \jsum x_j(1 - E[Z_{j,n} | \F_{n-1}]).
\end{align*}
Note that some of the variables in the right hand side of the last equation are random, but are almost surely constant given $ \F_{n-1} $. We thus need to show that, for $ x,y \in \realnonnegative^N, \lambda \in [0,1]$,
\[
	f_n\left(\lambda x + (1-\lambda)y \right) \leq \lambda f_n(x) + (1-\lambda)f_n(y).
\]
Note $ \sigma_i(x) $ is linear in $ x $. 
Moreover,
\begin{align*}
	&f_n\left(\lambda x + (1-\lambda)y \right) - \lambda f_n(x) - (1-\lambda)f_n(y) \\
	&= \frac{1}{N}\Nsum\frac{c_i}{d_i + \sigma_{i}(\lambda x + (1-\lambda)y)} - \frac{\lambda c_i}{d_i + \sigma_{i}(x)} - \frac{(1 - \lambda)c_i}{d_i + \sigma_{i}(y)} \\
	&= \frac{1}{N}\quad\mathclap{\Nsum}\enskip\,\frac{c_i\lambda(\lambda - 1)(\sigma_{i}(x) - \sigma_{i}(y))^2}{(d_i + \sigma_{i}(\lambda x + (1-\lambda)y))(d_i + \sigma_{i}(x))(d_i + \sigma_{i}(y))} \\
	&\leq 0,
\end{align*}
since $ \lambda - 1 \leq 0 $ and all other terms are nonnegative. As a result, $ E[\exposure | \F_{n-1}] $ is convex in the curing parameters $ (\Delta_{b,1}(n),\ldots,\Delta_{b,N}(n)) $ for all time. Lastly, the constraint set $ \left\{\{\Delta_{b,i}(n)\}_{i=1}^N \in \realnonnegative^N \ | \ \Nsum \Delta_{b,i}(n) = \budget \right\} $ is clearly a finite-dimensional simplex and hence convex and compact. \oprocend

\end{appendix}

\section{Acknowledgements}\label{sec:acknowledgements}

The authors wish to acknowledge the Centre for Advanced Computing at Queen's University, whose computing cluster allowed the simulations presented herein to be performed.

\bibliographystyle{ieeetr}
\bibliography{alias,Main-add,MH-add}
\end{document}